\documentclass[preprint,12pt,3p]{elsarticle}

\usepackage{amssymb}
\usepackage{amsmath,amsthm}
\usepackage{mathrsfs}
\usepackage{hyperref}
\usepackage{cleveref}
\usepackage[all,cmtip]{xy}
\usepackage{epsf, graphicx}
\usepackage{latexsym,amsfonts,amsbsy,amssymb}
\usepackage{geometry}
\usepackage{titletoc}
\usepackage{color,xcolor}
\usepackage{rotating}
\usepackage{algorithm}
\usepackage{algorithmic}
\usepackage{amscd}

\makeatletter

\@addtoreset{equation}{section} \makeatother
\newtheorem{theorem}{Theorem}[section]
\newtheorem{lemma}[theorem]{Lemma}

\newtheorem{void}[theorem]{}
\newtheorem{proposition}[theorem]{Proposition}

\def\Res{{\rm Res}}
\def\Def{{\rm Def}}
\def\Inf{{\rm Inf}}

\def\Aut{{\rm Aut}}

\def\End{{\rm End}}

\def\G{\mathcal{G}}

\def\O{\mathcal{O}}

\def\F{\mathcal{F}}
\def\G{\mathcal{G}}

\makeatletter
\def\ps@pprintTitle{%
\let\@oddhead\@empty
\let\@evenhead\@empty
\def\@oddfoot{\reset@font\hfil\thepage\hfil}
\let\@evenfoot\@oddfoot
}
\makeatother

\begin{document}

\begin{frontmatter}

\title{On stable equivalences with endopermutation source and K\"ulshammer--Puig classes}

\author{Xin Huang}


\begin{abstract}
We give a new proof, by using the terminology and notation in the textbook \cite{Lin18b}, to a result, due to Puig, stating that a stable equivalence of Morita type between two block algebras of finite groups induced by a bimodule with an endopermutation source preserves K\"ulshammer--Puig classes.
\end{abstract}


\end{frontmatter}


\section{Introduction}\label{s1}

In addition to defect groups and fusion systems, the third fundamental local invariant of a block is the family of {\it K\"{u}lshammer--Puig classes} of this block; we refer to \cite[Section 8.4]{Lin18b} or \ref{void:KP class} below for the definition. The number of conjugacy classes of weights in a block can be expressed in terms of a fusion system and the K\"ulshammer--Puig classes of that block; see \cite[Proposition 5.4]{Kessar} or \cite[Theorem 8.14.4]{Lin18b}.


Throughout this paper $p$ is a prime and $k$ is a field of characteristic $p$, which is large enough for finite groups considered below.


Let $G$ be a finite group, $b$ a block idempotent of $kG$, and $(D,e_D)$ a maximal $b$-Brauer pair of $kG$. For any subgroup $P$ of $D$, let $e_P$ be the unique block idempotent of $kC_G(P)$ such that $(P,e_P)\leq (D,e_D)$. (See \cite[Section 6.3]{Lin18b} for details on Brauer pairs). Let $\F=\F_{(D,e_D)}(G,b)$ be the fusion system of $kGb$ with respect to the choice of $(D,e_D)$; see \cite[8.5.1 and 8.5.2]{Lin18b}. Let $H$ be another finite group, $c$ a block idempotent of $kH$ and $(E,f_E)$ a maximal $c$-Brauer pair of $kH$. Similarly, for any subgroup $Q$ of $E$, let $f_Q$ be the unique block idempotent of $kC_H(Q)$ such that $(Q,f_Q)\leq (E,f_E)$, and let $\G:=\F_{(E,f_E)}(H,c)$. Assume that $M$ is an indecomposable
$(kGb, k Hc)$-bimodule, which is finitely generated projective as a left and right module, inducing a stable equivalence of Morita type between $k Gb$ and $k Hc$. That means that $M\otimes_{kHc}M^*\cong kGb\oplus W$ for some projective $kGb\otimes_k (kGb)^{\rm op}$-module $W$ and $M^*\otimes_{kGb} M\cong kHc\oplus W'$ for some projective $kHc\otimes_k (kHc)^{\rm op}$-module $W'$, where $M^*:={\rm Hom}_k(M,k)$. 
Assume further that as a $k(G\times H)$-module, $M$ has an endopermutation module (see e.g. \cite[Definition 7.3.1]{Lin18b}) as a source for some vertex. In \cite[7.6]{Puig1999}, Puig stated, using different terminology and notation, that various invariants of $kGb$ and $kHc$ are preserved under such a stable equivalence, including defect groups, fusion systems, and K\"ulshammer--Puig classes. The invariance properties for defect groups and fusion systems were proved by Linckelmann in \cite[Theorem 9.11.2]{Lin18b} using different terminology, notation and methods:

\begin{theorem}[cf. {{\cite[7.6.3]{Puig1999}}, see \cite[Theorem 9.11.2]{Lin18b}}]\label{Puig1}
Keep the notation and assumptions above. Then there is an isomorphism $\varphi : D\to E$ such that ${^\varphi{\F}}=$ $\mathcal{G}$, where ${}^\varphi\F$ is the fusion system on $E$ induced by $\F$ and the isomorphism $\varphi$. 


	
	
	
\end{theorem}


In this paper we give a new proof for the invariance of K\"ulshammer--Puig classes using the terminology and notation in \cite{Lin18b}. The main result we are going to prove is the following:

\begin{theorem}[cf. {\cite[7.6.5]{Puig1999}}]\label{theo:main result}
Keep the notation and assumptions above. Let $P$ be an $\F$-centric subgroup of $D$ - this is equivalent to require $Z(P)$ to be a defect group of the block algebra $kC_G(P)e_P$. Then $Q:=\varphi(P)$ is a $\G$-centric subgroup of $E$ by Theorem \ref{Puig1} and the preservation of centric subgroup. Let $\alpha_P\in H^2(\Aut_\F(P),k^\times)$ and $\beta_Q\in H^2(\Aut_\G(Q),k^\times)$ denote the K\"ulshammer--Puig classes of $(P,e_P)$ and $(Q,f_Q)$, respectively. Denote by $\eta:\Aut_\F(P)\xrightarrow{\sim}\Aut_\G(Q)$ the isomorphism induced by $\varphi$ in Theorem \ref{Puig1}. Then 
$$\beta_Q=\eta^*(\alpha_P)\in H^2(\Aut_\G(Q),k^\times),$$
where $\eta^*$ denotes the obvious map from $H^2(\Aut_\F(P),k^\times)$ to $H^2(\Aut_\G(Q),k^\times)$ induced by the isomorphism $\eta$.
\end{theorem}

We give the new proof of Theorem \ref{theo:main result} in Section \ref{section:Proof} after reviewing some necessary notation and lemmas in Section \ref{section:Preliminaries}. Our proof of Theorem \ref{theo:main result} can be viewed as another application of the main result Theorem 1.2 in \cite{H23}. The methods in this proof are also inspired by the proof of \cite[Theorem 1.5]{BP}.

\section{Preliminaries}\label{section:Preliminaries}

For a $k$-algebra $A$, we denote by $A^{\rm op}$ the opposite $k$-algebra of $A$. Unless specified otherwise, modules in the paper are left modules and are finite-dimensional over $k$. 
The following easy lemma will be used later.

\begin{lemma}\label{lem:Morita-bimodule}
	Let $A$ and $B$ be two $k$-algebras. Assume that an $A$-$B$-bimodule $M$ and a $B$-$A$-bimodule $N$ induce a Morita equivalence between $A$ and $B$. Then the $(A\otimes_k A^{\rm op})$-$(A\otimes_k B^{\rm op})$-bimodule $A\otimes_k N$ and the $(A\otimes_k B^{\rm op})$-$(A\otimes_k A^{\rm op})$-bimodule $A\otimes_k M$  induces a Morita equivalence between $A\otimes_k A^{\rm op}$ and $A\otimes_k B^{\rm op}$. 
\end{lemma}

\noindent{\it Proof.} This can be checked straightforward by definition.   $\hfill\square$

\begin{void}
	{\rm \textbf{Some usual convention and notation on modules and blocks.}
Let $G$ and $H$ be finite groups. We denote by $\Delta G$ the diagonal subgroup $\{(g,g)|g\in G\}$ of the direct product $G\times G$. A $k G$-$kH$-bimodule $M$ can be regarded as a $k(G\times H)$-module (and vice versa) via $(g,h)m=gmh^{-1}$, where $g\in G$, $h\in H$ and $m\in M$. If $M$ is indecomposable as a $k G$-$kH$-bimodule, then $M$ is indecomposable as a $k(G\times H)$-module, hence has a vertex (in $G\times H$) and a source.

By a {\it block} of the group algebra $k G$, we mean a primitive idempotent $b$ in the center of $k G$, and $k Gb$ is called a {\it block algebra} of $k G$. 
Denote by $-^\circ$ the $k$-algebra isomorphism $k H\cong(k H)^{\rm op}$ sending any $h\in H$ to $h^{-1}$. Let $b$ and $c$ be blocks of $k G$ and $k H$ respectively. Clearly $c^\circ$ is a block of $k H$. Then a $kGb$-$kHc$-bimodule $M$ can be regarded as a $k(G\times H)$-module belonging to the block $b\otimes c^\circ$ of $k(G\times H)$, and vice versa. Here, we identify $b\otimes c^\circ$ and its image under the $k$-algebra isomorphism $kG\otimes_k  kH\cong k(G\times H)$ sending $g\otimes h$ to $(g,h)$ for any $g\in G$ and $h\in H$. For any $b$-Brauer pair $(P,e_P)$ of $k Gb$, we set $N_G(P,e_P):=\{g\in N_G(P)~|~ge_Pg^{-1}=e_P\}$. By \cite[Theorem 6.2.6 (iii)]{Lin18b}, $e_P$ remains a block of $kN_G(P,e_P)$.
}
\end{void}

\begin{void}[see {\cite{Bouc10book}} or {\cite[2.1, 6.1(a)]{BP}}]\label{void:Subgroups of direct product groups}
{\rm \textbf{Subgroups of direct product groups.} Let $G$, $H$, $K$ be finite groups, $X$ a subgroup of $G\times H$, and $Y$ a subgroup of $H\times K$. We denote by $p_1:G\times H\to G$ and $p_2:G\times H\to H$ the canonical projections. Let 
	$$k_1(X):= \{g\in G~|~(g,1)\in X\}~~{\rm and}~~k_2(X):=\{h\in H~|~(1,h)\in X\}.$$
Then one obtains normal subgroups $k_i(X)$ of $p_i(X)$ and canonical isomorphisms 
$X/(k_1(X)\times k_2(X))\to p_i(X)/k_i(X)$, for $i\in \{1,2\}$, induced by the projection maps $p_i$. The resulting isomorphism $$\eta_X:p_2(X)/k_2(X)\xrightarrow{\sim}p_1(X)/k_1(X)$$ satisfies $\eta_X(hk_2(X))=gk_1(X)$ if and only if $(g,h)\in X$. Here $(g,h)\in p_1(X)\times p_2(X)$.
Let the notation $X^\vee$ denote the subgroup
$\{(h,g)~|~(g,h)\in X\}$ of $H\times G$. The composition of $X$ and $Y$ is defined as
$$X*Y:=\{(g,k)\in G\times K~|~\exists h\in H, (g,h)\in X, (h,k)\in Y\}.$$
This is a subgroup of $G\times K$. Let $M$ be a $kX$-module and $N$ a $kY$-module. Since $k_1(X)\times k_2(X)\leq X$, $M$ can be viewed as a $k(k_1(X))$-$k(k_2(X))$-bimodule via restriction. Similarly, $N$ can be viewed as a $k(k_1(Y))$-$k(k_2(Y))$-bimodule, and $M\otimes_{k(k_2(X)\cap k_1(Y))}N$ is a $k(k_1(X))$-$k(k_2(Y))$-bimodule. Note that $k_1(X)\times k_2(Y)\leq X*Y$ and that this bimodule can be extended to a $k(X*Y)$-module such that for any $(g,k)\in X*Y$, $m\in M$ and $n\in N$,
$$(g,k)(m\otimes n)=(g,h)m\otimes (h,k)n,$$
where $h\in H$ is any element such that $(g,h)\in X$ and $(h,k)\in Y$. This $k(X*Y)$-module is denoted by $M\mathop{\otimes}\limits_{kH}^{X,Y}N$.
}
\end{void}

\begin{void}\label{void:KP class}
{\rm \textbf{Schur classes and K\"{u}lshammer--Puig classes.} Let $G$ be a finite group, $N$ a normal subgroup of $G$, and $S$ a $G$-stable simple $kN$-module. Since we assumed that $k$ is large enough, we have $\End_{kN}(S)=k$. These data gives rise to a canonical cohomology class $\alpha\in H^2(G/N,k^\times)$ which is called the {\it Schur class} of $S$ with respect to $N\unlhd G$; see e.g. \cite[Theorem 5.3.12 (i), (ii)]{Lin18a} for the construction of $\alpha$. It is well-known and easy to check that the Schur class of the $k$-dual $S^*$ of $S$ with respect to $N\unlhd G$ is $\alpha^{-1}$.

Let $G$ be a finite group, $b$ a block of $kG$, and $(D,e_D)$ a maximal $b$-Brauer pair. Denote by $\F$ the fusion system of the block $b$ on it defect group $D$ determined by the choice of $(D,e_D)$. For any subgroup $P$ of $D$ denote by $e_P$ the unique block of $kC_G(P)$ satisfying $(P,e_P)\leq (D,e_D)$. For any $\F$-centric subgroup $P$ of $D$, since $Z(P)$ is a defect group of $kC_G(P)e_P$ (see e.g. \cite[Proposition 8.5.3 )(iii)]{Lin18b}), $kC_G(P)e_P$ has a unique simple module $S$; see e.g. \cite[Proposition 6.6.5]{Lin18b}. Hence $S$ is $N_G(P,e_P)$-stable. Considering the Schur class of $S$ with respect to $C_G(P)\unlhd N_G(P,e_P)$, and noting that $\Aut_\F(P)\cong N_G(P,e_P)/C_G(P)$, we can identify this Schur class with an element $\alpha_P\in H^2(\Aut_\F(P),k^\times)$. The class $\alpha_P$ is called the {\it K\"ulshammer--Puig classes} of $(P,e_P)$.
}
\end{void}

The following lemma is very useful in the proof of Theorem \ref{theo:main result}.

\begin{lemma}[{\cite[Lemma 11.7]{BP}}]\label{lemma:useful lemma from BP}
Let $G$ and $H$ be finite groups, and let $Y$ be a subgroup of $G\times H$ such that $p_1(Y)=G$ and $p_2(Y)=H$. Let $\eta_Y: H/k_2(Y)\xrightarrow{\sim} G/k_1(Y)$ be the isomorphism induced by $Y$ (see \ref{void:Subgroups of direct product groups}). Suppose that $S_1$ is a $G$-stable simple $k(k_1(Y))$-module and $S_2$ is an $H$-stable simple $k(k_2(Y))$-module. Denote by $\alpha\in H^2(G/k_1(Y),k^\times)$ and $\beta\in H^2(H/k_2(Y))$ their respective Schur classes. Suppose further that there exists a $kY$-module W such that $\Res_{k_1(Y)\times k_2(Y)}^Y(W)\cong S_1\otimes_k S_2^*$. Then 
$$\beta=\eta_Y^*(\alpha)\in H^2(H/k_2(Y)),$$
where $\eta_Y^*$ denotes the obvious map from $H^2(G/k_1(Y),k^\times)$ to $H^2(H/k_2(Y),k^\times)$ induced by the isomorphism $\eta_Y$.
\end{lemma}

\begin{void}
{\rm \textbf{Inflation and deflation functors.} Let $G$ be a finite group and $N$ a normal subgroup of $G$. The canonical homomorphism $\pi:G\to G/N$ induces two functors: The {\it inflation functor} $\Inf_{G/N}^G$ from the category of $k(G/N)$-modules to the category of $kG$-modules which is given by restriction along $\pi$; and the {\it deflation functor} $\Def_{G/N}^G$ in the other direction, which is just the tensor functor $k(G/N)\otimes_{kG}-$, where $k(G/N)$ is viewed as a $k(G/N)$-$kG$-bimodule using $\pi$ for the right module structure. 
	

}		
\end{void}

\section{Proof of Theorem \ref{theo:main result}}\label{section:Proof}

In this section we give our proof to Theorem \ref{theo:main result}. We first fix some notation.

\begin{void}\label{void:notation-slashmodules}
{\rm \textbf{Notion.}
Keep the notation and assumptions of Theorem \ref{Puig1}. We may identify $D$ and $E$ via the isomorphism $\varphi$. So now $D=E$, $\varphi$ is the identity automorphism of $D$ and $\F=\mathcal{G}$. (This simplifies notation, but one could as well go on without making those identifications.) Consider the block $b\otimes c^\circ$ of $k(G\times H)$. Then $(D\times D,e_D\otimes f_D^\circ)$ is a maximal $(b\otimes c^\circ)$-Brauer pair of $k(G\times H)$; see \cite[Proposition 4.9]{Lin08}. 
Again by \cite[Proposition 4.9]{Lin08}, the fusion system of $k(G\times H)(b\otimes c^\circ)$ with respect to the maximal $(b\otimes c^\circ)$-Brauer pair $(D\times D,e_D\otimes f_D^\circ)$ is $\F\times \G$; see \cite[page 1282]{Lin08} for the notation $\F\times \G$. The $k(G\times H)(b\otimes c^\circ)$-module $M$ is a Brauer-friendly module in the sense of \cite[Definition 8]{Biland}; see \cite[Proposition 2.6]{H23}.  For any subgroup $P$ of $D$, consider the $(b\otimes c^\circ)$-Brauer pair $(\Delta P,e_P\otimes f_P^\circ)$ of $k(G\times H)$. Let $M_P$ be a $(\Delta P,e_P\otimes f_P^\circ)$-slashed module attached to $M$ over the  group $N_{G\times H}(\Delta P, e_P\otimes f_P^\circ)$; see e.g. \cite[Lemma 2.8]{H23} for its definition. Note that $M_P$ is a $kN_{G\times H}(\Delta P, e_P\otimes f_P^\circ)$-module on which $\Delta P$ acts trivially. By \cite[Theorem 1.2 (a)]{H23}, $\Res_{C_G(P)\times C_G(P)}^{kN_{G\times H}(\Delta P,e_P\otimes f_P^\circ)}(M_P)$ induces a Morita equivalence between the block algebras $kC_G(P)e_P$ and $kC_H(P)f_P$.

}
\end{void}

\begin{proposition}\label{prop:extensions}
Keep the notation of \ref{void:notation-slashmodules}. Let $P$ be a subgroup of $D$ and suppose that $P$ is centric in $\F$, or equivalently, $Z(P)$ is a defect group of $kC_G(P)e_P$. Set 
$$I:=N_G(P,e_P),~~J:=N_H(P,f_P),~~Y:=N_{G\times H}(\Delta P,e_P\otimes f_P^\circ),~~C:=C_G(P)\times C_H(P),$$
$$Z:=Z(P)\times Z(P),~~{\rm and}~~Z':=\Delta Z(P).$$
The following hold.

\begin{enumerate} [\rm (a)]
	\item  ${\rm Res}_C^Y(M_P)$ is an indecomposable $kC$-module having $Z'$ as a vertex and a trivial $kZ'$-source.
	
	\item The $kC$-module $S:={\rm Inf}_{C/Z}^{C}{\rm Def}_{C/Z}^{C}{\rm Res}_C^Y(M_P)$ is the unique simple $kC(e_P\otimes f_P^\circ)$-module, up to isomorphism.
	
	\item There exists a simple $kY(e_P\otimes f_P^\circ)$-module $W$ such that ${\rm Res}_C^Y(W)=S$.
	
\end{enumerate}

\end{proposition}	

\noindent{\it Proof.} Since $b$ and $c$ have a common fusion system $\F=\G$, we see that $Z(P)$ is also a defect group of $\O C_H(P)f_P$. It follows that $Z$ is a defect group of $kC(e_P\otimes f_P^\circ)$; see e.g. \cite[Proposition 4.9 (i)]{Lin08}. By the definition of $M_P$, $\Delta P$ acts trivially on $M_P$, and hence $Z'$ acts trivially on $M_P$. So $Z'$ is contained in a vertex of ${\rm Res}_C^Y(M_P)$; see e.g. \cite[Theorem 5.6.9]{Lin18a}. Since ${\rm Res}_C^Y(M_P)$ has an endopermutation module as a source (see \cite[Theorem 1.2 (a)]{H23}) and induces a Morita equivalence, by \cite[Theorem 9.11.2 (i)]{Lin18b}, the order of any vertex of ${\rm Res}_C^Y(M_P)$ equals to the order of any defect group of $kC_G(P)e_P$. This implies that $Z'$ is a vertex of ${\rm Res}_C^Y(M_P)$. Since ${\rm Res}_{Z'}^{Y}(M_P)$ is a direct sum of trivial $kZ'$-modules, if $L$ is a $kZ'$-source of ${\rm Res}_C^Y(M_P)$, the trivial $kZ'$-module $k$ is isomorphic to a direct summand of ${\rm Res}_{Z'}^C{\rm Ind}_{Z'}^C(L)$. 
Using Mackey's formula and noting that $Z'$ is in the center of $C$, we easily deduce that $L$ is a trivial $kZ'$-module, whence (a).
Statement (b) follows from (a) and \cite[Proposition 4.8 (b)]{BP}.

It remains to prove (c). Recall that $Z'$ acts trivially on $M_P$. Using $Z=Z'\cdot(Z(P)\times \{1\})$ and \cite[Proposition 6.7 (d)]{BP}, we obtain
$$\Inf_{Y/Z}^Y\Def_{Y/Z}^Y(M_P)={\rm Inf}_{Y/(Z(P)\times\{1\})}^Y\Def_{Y/(Z(P)\times\{1\})}^Y(M_P)\cong k(C_G(P)/Z(P))\bar{e}_P\mathop{\otimes}\limits_{kG}^{Y*Y^\vee,Y}M_P,$$
where $\bar{e}_P$ denotes the image of $e_P$ under the canonical map $kC_G(P)\to k(C_G(P)/Z(P))$. By a straightforward calculation (or by using \cite[Lemma 2.2 (a)]{BP}), we see that $Y*Y^\vee=(\Delta I)\cdot (C_G(P)\times C_G(P))$.  


Denote by $W$ the $kY$-module $k(C_G(P)/Z(P))\bar{e}_P\mathop{\otimes}\limits_{kG}^{Y*Y^\vee,Y}M_P$. By definition, 
$$\Res_C^Y(W)=k(C_G(P)/Z(P))\bar{e}_P\otimes_{kC_G(P)} \Res_C^Y(M_P).$$ We claim that $W$ is simple. It suffices to show that ${\rm Res}_C^Y(W)$ is simple. Indeed, since $kC_G(P)e_P$ has a central defect group $Z(P)$, $k(C_G(P)/Z(P))\bar{e}_P$ is a block of defect zero; see e.g. \cite[Proposition 6.6.5]{Lin18b}. Thus $k(C_G(P)/Z(P))\bar{e}_P$ is a simple $k(C_G(P)/Z(P))\bar{e}_P$-$k(C_G(P)/Z(P))\bar{e}_P$-bimodule and via inflation a simple $k(C_G(P)\times C_G(P))(e_P\otimes e_P^\circ)$-module. Since ${\rm Res}_C^Y(M_P)$ induces a Morita equivalence between $kC_G(P)e_P$ and $kC_H(P)f_P$, the tensor functor $-\otimes_{kC_G(P)}{\rm Res}_C^Y(M_P)$ induces an equivalence between the category of $kC_G(P)e_P$-$kC_G(P)e_P$-bimodules and the category of $kC_G(P)e_P$-$kC_H(P)f_P$-bimodules; see Lemma \ref{lem:Morita-bimodule}. From this, we see that $\Res_C^Y(W)$ is simple, as claimed. Since clearly
$$\Res_C^Y\Inf_{Y/Z}^Y\Def_{Y/Z}^Y(M_P)\cong \Inf_{C/Z}^C\Def_{C/Z}^{C}\Res_C^Y(M_P),$$
now we have $\Res_C^Y(W)\cong\Res_C^Y\Inf_{Y/Z}^Y\Def_{Y/Z}^Y(M_P)\cong \Inf_{C/Z}^C\Def_{C/Z}^{C}\Res_C^Y(M_P)=S$. This completes the proof.  $\hfill\square$

\medskip\noindent{\it Proof of Theorem \ref{theo:main result}.} Since $P$ is an $\F$-centric subgroup of $D$, we are in the context of Proposition \ref{prop:extensions}, and hence we can use the notation there. Let $S_1$ be a simple $kC_G(P)e_P$-module and $S_2$ a simple $kC_H(P)f_P$-module. Then by \ref{void:KP class}, $\alpha_P$ is the Schur class of $S_1$ with respect to $C_G(P)\unlhd N_G(P,e_P)$ and $\beta_P:=\beta_Q$ is the Schur class of $S_2$ with respect to $C_H(P)\unlhd N_H(P,f_P)$. Consider the inclusion $Y\leq I\times J$, we have $p_1(Y)=I$,  $p_2(Y)=J$, $k_1(Y)=C_G(P)$ and $k_2(Y)=C_H(P)$. By an easy calculation we see that the resulting isomorphism $\eta_Y:p_1(Y)/k_1(Y) \xrightarrow{\sim} p_2(Y)/k_2(Y)$ (see \ref{void:Subgroups of direct product groups}) is exactly the isomorphism $\Aut_\F(P)\xrightarrow{\sim} \Aut_\G(P)$. By Lemma \ref{lemma:useful lemma from BP}, it suffices to show that the simple $kC(e_P\otimes f_P^\circ)$-module $S_1\otimes_k S_2^*$ can be extended to $kY$. But this has been proven in Proposition \ref{prop:extensions}.   $\hfill\square$



\begin{thebibliography}{1}
\expandafter\ifx\csname url\endcsname\relax
  \def\url#1{\texttt{#1}}\fi
\expandafter\ifx\csname urlprefix\endcsname\relax\def\urlprefix{URL }\fi
\expandafter\ifx\csname href\endcsname\relax
  \def\href#1#2{#2} \def\path#1{#1}\fi






\bibitem{Biland}
E. Biland, Brauer-friendly modules and slash functors, J. Pure Appl. Algebra {\bf218} (2014) 2319--2336.






\bibitem{BP}
R. Boltje, P. Perepelitsky, $p$-Permutation equivalences between blocks of group algebras, J. Algebra {\bf 664} (2025) 815--887.


\bibitem{Bouc10book}
S. Bouc, Biset Functors for Finite Groups, Lecture Notes in Mathematics, vol. 1990, Springer-Verlag, Berlin, 2010.










\bibitem{H23}
X. Huang, Stable equivalences with endopermutation source, slash functors, and functorial equivalences, J. Algebra {\bf 650} (2024) 335--353. 


\bibitem{Kessar}
R. Kessar, Introduction to block theory, Group Representation Theory, EPFL Press, Lausanne (2007) 47--77.



\bibitem{Lin08}
M. Linckelmann, Trivial source bimodule rings for blocks and $p$-permutation equivalences, Trans. Amer. Math. Soc. {\bf361} (2009) 1279--1316.





\bibitem{Lin18a}
M. Linckelmann, The Block Theory of Finite Group Algebras I, London Math. Soc. Student Texts, vol. {\bf91}, Cambridge University Press, 2018.

\bibitem{Lin18b}
M. Linckelmann, The Block Theory of Finite Group Algebras II, London Math. Soc. Student Texts, vol. {\bf92}, Cambridge University Press, 2018.




\bibitem{Puig1999}
L. Puig, On the Local Structure of Morita and Richard Equivalences Between Brauer Blocks, Progress in Math., vol. {\bf 178}, Birk\"{a}user Verlag, Basel, 1999.









\end{thebibliography}
\end{document}